\documentclass[10pt]{article}
\usepackage{amsmath}
\usepackage{amsfonts}
\usepackage{amssymb}
\usepackage{amscd}
\usepackage{amsthm}
\usepackage{latexsym}
\usepackage{amsbsy}
\newtheorem{lem}{Lemma}[section]
\newtheorem{prop}{Proposition}[section]
\newtheorem{theorem}{Theorem}[section]

\newtheorem{definition}{Definition}[section]

\newcommand{\sw}[1]{^{(#1)}}

\newcommand{\ot}{\otimes}
\begin{document}
\title{A Note on Cyclic Duality and Hopf Algebras}
\author {\hspace{-1cm}M. Khalkhali~~~~~~~~~~~~~~~~~~~~~B. Rangipour\thanks{PIMS postdoctoral  fellow}
\\ \footnotesize\text{\hspace{-1cm}masoud@uwo.ca~~~~~~~~~~~~~~~~~~~~~~~~~~~~~~~~~~~bahram@uvic.ca}
\\\footnotesize{Department of Mathematics~~~~~~~~~ Department of Mathematics and Statistics}
\\\footnotesize{\hspace{-1.8cm}University of Western Ontario~~~~~~~~~~~~~~~~~~~~~~University of Victoria} }
\maketitle

\begin{abstract}
We show that various cyclic and cocyclic modules attached to Hopf
algebras and Hopf modules are related to each other via Connes' duality isomorphism
for the cyclic category.
 \end{abstract}
 \section{Introduction}

A remarkable property of Connes' cyclic category $\Lambda$ is its self duality in the sense that there is a
natural isomorphism  between $\Lambda$ and
its opposite  category  $\Lambda^{op}$  \cite{ac88}. Roughly  speaking, the duality functor $\Lambda^{op} \longrightarrow
\Lambda $ acts as identity on objects
of $\Lambda$ and exchanges face and degeneracy operators while sending the cyclic operator to its inverse
(see Section 2 for
a precise formulation). Thus to a cyclic (resp. cocyclic) module one can
associate a cocyclic (resp. cyclic) module by applying Connes' duality
isomorphism. This notion of duality, called cyclic duality in this paper,  should not be confused with the
duality obtained by applying the $Hom$ functor, since it is of a very
different nature. For example, it is easy to see that the cyclic dual of the cyclic
(resp. cocyclic)
module of an algebra (resp. coalgebra) is homologically trivial (Lemma
2.3
below).

The goal of this article is to show that the situation is much more
interesting for Hopf algebras by showing that various, non-trivial,  cyclic and cocyclic
modules attached to
Hopf algebras and Hopf modules are cyclic duals of each other.  Recently we have seen  a
proliferation of cyclic and cocyclic modules attached
to Hopf algebras \cite{cr, ta, kr, kr5, hkrs, js}, extending the pioneering work of Connes and Moscovici
\cite{achm98, achm00}. Recall that in \cite{hkrs} a cocyclic
module $C^*_H (H, M)$ and a cyclic module $C^H_* (H, M)$ is defined for any Hopf algebra $H$ and an
stable anti-Yetter-Drinfeld $H$-module $M$ (see Section 3 for definitions). For
$M=k$, the ground field, these modules reduce to the Connes-Moscovici
cocyclic module of an Hopf algebra endowed with a modular pair in
involution \cite{achm98, achm00} and  to the cyclic module introduced in \cite{kr, ta},
respectively. It is known that the associated cyclic homology and
cohomology theories are non-trivial.
In this paper we show that  $C^H_* (H, M)$ is in fact
isomorphic, via a non-trivial map, to the cyclic dual of $C^*_H (H, M)$ (Theorem 3.1).

The cyclic dual of $C^*_H (H, M)$ appears naturally in the study of the
relative cyclic homology of Hopf-Galois extensions \cite{js}. Thus our result
 shows that this theory is a special case of the invariant cyclic
 homology defined in \cite{hkrs} for general coefficients and in \cite{kr5} for a restricted class of
 coefficients.

 We would like to thank Nigel Higson for valuable discussions on the
 subject of this paper. Our thanks go also to our collaborators in
 \cite{hkrs}, Piotr. M. Hajac and Yorck Sommerh\"auser.

\section{Duality for cyclic modules}

Let $k$ be a commutative ring with identity and let $k-mod$ denote the category of $k$-modules.  Recall that a cyclic
$k$-module (or a cyclic module for short)  is
a contravariant functor $\Lambda \rightarrow k-mod$, where $\Lambda $ denotes
 Connes' category \cite{ac88} (cf. also \cite{ac, ld}). Equivalently, a
cyclic module is given by
a sequence $X_n, n \geq 0$, of $k$-modules  and $k$-linear maps  called face,
 degeneracy and cyclic operators
 $$\delta_i: X_{n} \rightarrow X_{n-1},\quad  \sigma_i : X_n \rightarrow X_{n+1},
 \quad  \tau_n : X_n \rightarrow X_n \qquad  0\le i\le n,$$
 such that $(X,\delta_i,\sigma_i)$ is a simplicial module  and the following extra
  relations are satisfied:
\begin{eqnarray*}
\delta_i \tau_n&=& \tau_{n-1}\delta_{i-1}\hspace{43 pt} 1\le i\le
n,\\
\delta_0 \tau_n&=& \delta_{n},
\\ \sigma_i \tau_n&=& \tau_{n-1}\sigma _{i-1} \hspace{43 pt} 1\le i\le n,\\
\sigma_0 \tau_n&=& \tau_n^2\sigma_n, \\
\tau_n^{(n+1)} &=& \mbox{id} _n .
\end{eqnarray*}
When all relations, except possibly the last one,  are satisfied we say that
we have a $paracyclic$ module.

A $cocyclic$ module is a functor $\Lambda \longrightarrow k-mod$.
Equivalently, a cocyclic module is given by a sequence of
$k$-modules $X^n$ and $k$-linear maps  called coface, codegeneracy, and
cyclic operators:
$$d_i: X^n \rightarrow X^{n+1},\quad  s_i : X^n \rightarrow X^{n-1},
 \quad  t_n : X^n \rightarrow X^n \qquad  0\le i\le n,$$
such that $(X, d_i, s_i)$ is a cosimplicial module and the following
extra relations are satisfied:
\begin{eqnarray*}
t_{n+1}d_i &=& d_{i-1} t_n \hspace{43 pt} 1\le i\le
n,\\
t_{n+1} d_0 &=& d_{n},
\\ t_{n-1} s_i &=& s_{i-1} t_n \hspace{43 pt} 1\le i\le n,\\
t_{n-1} s_0 &=& s_n t_n^2, \\
t_n^{(n+1)} &=& \mbox{id} _n .
\end{eqnarray*}
When all  relations, except possibly the last one, are satisfied we say that
we have a $paracocyclic$ module.

Let $X=(X^n, d_i, s_i, t_n)$ be a paracocyclic (resp. cocyclic)
module where we assume that $t_n$ is invertible for all $n\geq 0$. We
denote its  cyclic dual by $\widehat{X}$. It is defined as follows \cite{ac88} . We put
$\widehat{X}_n=X^n$, and
\begin{eqnarray*}
  \delta_i &=&s_{i-1}: \widehat{X}_n \longrightarrow \widehat{X}_{n-1} , \hspace{43 pt} 1\leq i \leq n,\\
\delta_0 &=&s_{n-1}t_n,\\
\sigma_i &=& d_i : \widehat{X}_n \longrightarrow \widehat{X}_{n+1},\\
\tau_n &=& t_n^{-1}.
\end{eqnarray*}

The following two lemmas  are  proved in \cite{ac88}. Both can be
checked directly.

\begin{lem}
$\widehat{X}=(\widehat{X}_n, \delta_i, \sigma_i, \tau_n)$ is a paracyclic
module. If $X$ is a cocyclic module, then $\widehat{X}$ is a  cyclic module.
\end{lem}

Conversely, one can obtain from a paracyclic (resp. cyclic) module a
paracocyclic (resp. cocyclic) module as follows. Let $X=(X_n, \delta_i,
\sigma_i, \tau_n)$
 be a paracyclic module. We denote the cyclic dual of $X$ by
$\check{X}$ where $\check{X}^n=X_n$ and its coface, codegeneracy and cyclic operators
are defined by
\begin{eqnarray*}
d_i &=&\sigma_{i-1}: \check{X}^n \longrightarrow \check{X}^{n+1}  \hspace{43 pt}  1\leq i\leq n,\\
d_0 &=&\tau_n \sigma_n,\\
s_i &=&\delta_i: \check{X}^n  \longrightarrow \check{X}^{n-1} \hspace{43 pt}  0\leq i \leq n-1,\\
t_n &=&\tau_n^{-1}.
\end{eqnarray*}

\begin{lem} $\check{X}=(\check{X}^n, d_i, s_i, t_n)$ is a paracocyclic module. If $X$ is
a cyclic module, then $\check{X}$ is a  cocyclic module.
\end{lem}

We give a few examples of cyclic and cocyclic modules that will be used in this paper. The
 first example is the most fundamental example which motivated the whole theory \cite{ac88}.
\begin{itemize}
\item[1.]  Let $A$ be an algebra. The cyclic module $A^ \natural$
is defined by $A_n^\natural=A^{\otimes(n+1)} , n\geq 0$, with  face,  degeneracy
and cyclic operators  defined by
\begin{eqnarray*}
\delta_i(a_0 \otimes a_1\otimes \dots \otimes a_n)&=&a_0 \otimes \dots \otimes
a_{i}a_{i+1}\otimes \dots \otimes a_n,\\
\delta_n(a_0 \otimes a_1\otimes\dots \otimes a_n)&=&a_na_0 \otimes a_1 \otimes
\dots  \otimes a_{n-1},\\
\sigma_i(a_0 \otimes a_1\otimes\dots \otimes a_n)&=&a_0 \otimes\dots \otimes
a_i \otimes  1 \otimes \dots\otimes
 a_n,\\
\tau_n(a_0 \otimes a_1\otimes\dots \otimes a_n)&=& a_n \otimes a_0\dots \otimes a_{n-1}.
\end{eqnarray*}

\item[2.] Let $C$ be a coalgebra. The cocyclic module $C_\natural$ is
defined by $C_\natural ^n=C^{\otimes (n+1)},\;\; n\geq 0$, with coface, codegeneracy
 and cyclic operators given by:
\begin{eqnarray*}
d_i(c_0 \otimes c_1 \otimes \dots \otimes c_n)&=& c_0 \otimes \dots
\otimes  c_i^{(1)}\otimes c_i^{(2)}\otimes c_n ~~~ 0 \leq i \leq n,\\
d_{n+1}(c_0 \otimes c_1 \otimes \dots \otimes c_n)&=& c_0^{(2)}\otimes c_1
\otimes \dots \otimes c_n  \otimes c_0^{(1)},\\
s_i(c_0 \otimes c_1 \otimes \dots \otimes c_n)&=&c_0 \otimes \dots
c_i \otimes \varepsilon(c_{i+1})\otimes
\dots\otimes c_n ~~~0\leq i \leq n-1,\\
t_n(c_0 \otimes c_1 \otimes \dots \otimes c_n)  &=& c_1 \otimes c_2 \otimes
\dots \otimes c_n \otimes c_0,
\end{eqnarray*}
where as usual $\Delta(c)= c^{(1)}\otimes c^{(2)} $ denotes the
coproduct of $C$ (Sweedler's notation), and $\varepsilon $ denotes the
counit of $C$.

\item[3.] Let $H$ be a Hopf algebra, $\delta : H \longrightarrow k$ an
algebra map (a character) and $\sigma \in H$ a grouplike element.
Following
  \cite{achm98, achm00}, we say $(\delta,\sigma)$
is a {\it modular pair } if $\delta (\sigma) =1$ and a {\it  modular pair in involution }
if,  in addition,  $(\sigma^{-1}\widetilde{S})^2=id_{H} $ where  the {\it
twisted antipode }  $\tilde{S}$ is defined by
$$\widetilde{S}(h)= \delta(h^{(1)}) S(h^{(2)}).$$
 Given $H$ and a modular pair in involution $(\delta,\sigma )$ as above,
Connes and Moscovici define a cocyclic module  $H_{(\delta,\sigma)}^\natural$
 as follows.
Let   $H_{(\delta,\sigma)}^{\natural,0}=k$ and
$H_{(\delta,\sigma)}^{\natural,n}=H^{\otimes n}$ , $n\geq 1$.
The coface, codegeneracy and cyclic operators $d_i$, $s_i$, $t_n$ are
defined by
\begin{eqnarray*}
d_0(h_1 \otimes \dots \otimes h_n)&=& 1_\mathcal{H} \otimes h_1 \otimes\dots
 \otimes h_n \\
d_i(h_1 \otimes\dots \otimes h_n)&=& h_1 \otimes \dots \otimes\Delta (h_i)
\otimes\dots \otimes h_n
  \;\;\text{for}\;\;1\leq i \leq n \\
d _{n+1}(h_1 \otimes\dots \otimes h_n)&=& h_1 \otimes\dots \otimes h_n
\otimes \sigma \\
s_i(h_1 \otimes\dots \otimes h_n)&=& h_1 \otimes\dots \otimes\epsilon
(h_{i+1})\otimes\dots \otimes h_n
 \;\;\text{for}\;\;0 \leq i \leq n \\
t_n(h_1 \otimes\dots \otimes h_n )&=&\Delta^{n-1}\widetilde{S}(h_1)\cdot(h_2 \otimes\dots \otimes h_n \otimes 
\sigma).
\end{eqnarray*}

The cyclic cohomology of
  this cocyclic module is, by definition, the cyclic cohomology of the Hopf algebra $H$ with
  respect to $(\delta , \sigma)$.

 \item[4.]

In \cite{kr} and, independently,  \cite{ta}, a  cyclic module is
 associated to any Hopf algebra $H$ endowed with a modular pair in
 involution denoted by $\widetilde{H}^{(\delta,\sigma)}_{\natural}$.
  We have
  $\widetilde{H}^{(\delta,\sigma)}_{\natural , n} = H^
 {\otimes n} $, for $n> 0$  and $ \widetilde{H}_{\natural , 0}^{( \delta ,\sigma)}= k $.
 Its face,  degeneracy, and cyclic operators
 are as follows: \\
\begin{eqnarray*}
 {\delta}_0 (h_1 \otimes h_2 \otimes  \dots  \otimes h_n) &=&\epsilon (h_1)h_2 \otimes h_3
  \otimes \dots  \otimes h_n \\
 {\delta}_i (h_1 \otimes h_2 \otimes
 \dots  \otimes h_n) &=& h_1\otimes h_2 \otimes \dots \otimes
h_i h_{i+1}\otimes \dots  \otimes h_n   \\
{\delta}_n(h_1 \otimes h_2 \otimes  \dots  \otimes h_n) &=& \delta
(h_n)h_1\otimes
 h_2 \otimes \dots  \otimes h_{n-1} \\
\hspace{2cm}{\sigma}_i(h_1 \otimes h_2 \otimes \dots  \otimes h_n)
&=&h_1 \otimes h_2 \dots \otimes h_i \otimes 1 \otimes
 h_{i+1} \dots  \otimes h_n     \\
\tau _n(  h_1 \otimes h_2 \otimes  \dots  \otimes h_n ) &=&
 \delta(h_n^{(2)})\sigma S ( h_1^{(1)} h_2^{(1)} \dots h_{n-1}^{(1)}h_n^{(1)} )
\otimes \\
& &h_1^{(2)} \otimes \dots \otimes h_{n-1}^{(2)} .
\end{eqnarray*}

\end{itemize}

It is natural to ask what is the relation between the Hochschild and
cyclic homology groups of a cyclic module $X$ and the Hochschild and
cyclic cohomology groups of its dual cocyclic module $\check{X}$. The following
simple lemma answers this question for algebras and coalgebras. In the next section we answer
this question for Hopf algebras.

\begin{lem}
Let $A$ be a unital algebra over a field $k$. Then  the Hochschild cohomology of the cocyclic
 module $\check{A^{\natural}}$ is trivial in positive dimensions. Similarly, if $C$ is a
 coalgebra over $k$, then the Hochschild homology
 groups of the dual cyclic module $\widehat{C_{\natural}}$ are trivial in positive dimensions.
\end{lem}
\begin{proof}
Let $\phi$ be a linear functional on $A$ such that $\phi(1)=1$. One can easily check
that the following defines a contracting homotopy  for the  Hochschild
complex of $\check{A^\natural}$:
\begin{align*}
&h:A^{\ot(n+1)}\longrightarrow A^{\ot n},
&h(a_0\ot a_1\ot\cdots\ot a_n)=\phi(a_0)a_1\ot a_2\ot \cdots\ot a_n.
\end{align*}

In the coalgebra case  let $c$ be an element of $C$ such that
$\epsilon(c)=1$. We  define  a contracting homotopy  for the
Hochschild complex of $\hat{C_\natural}$ as follows:
\begin{align*}
&s:C^{\ot n}\longrightarrow C^{\ot(n+1)},
&&s(c_0\ot c_1\ot\dots\ot c_{n-1})=c\ot c_0\ot c_1\ot\dots \ot c_{n-1}.
\end{align*}

\end{proof}

\section{Cyclic duality and Hopf algebras}

To define a cyclic (co)homology theory with
coefficients for Hopf algebras, the module of coefficients must be of a very special type. In
\cite{hkrs} the most general allowable Hopf modules of this type are identified and called
stable anti-Yetter-Drinfeld (SAYD) modules. One dimensional SAYD modules
correspond exactly to modular pairs in involution. An intermediate case
are matched and comatched pairs of \cite{kr5}. In this section we first recall the notion of an SAYD module over a
Hopf algebra
and their associated cyclic and cocyclic modules from \cite{hkrs}. We
then prove that, quite unexpectedly, these modules are, up to isomorphism,  cyclic duals of
each other.

If $M$ is a left $H$-comodule we write $_M\Delta (m)=m^{(-1)} \otimes
m^{(0)}$ to denote its coaction $_M\Delta : M \rightarrow H \otimes M$
(Sweedler's notation). Similarly if $M$ is a right $H$-comodule, we
write $\Delta_M (m)= m^{(0)} \otimes m^{(1)}$ to denote its coaction
$\Delta_M : M \rightarrow M\otimes H$.

 \begin{definition}
Let $H$ be a Hopf algebra with a bijective antipode $S$,
 and $M$ a module and comodule
 over $H$. We call $M$ an anti-Yetter-Drinfeld module  if
 the action and coaction are compatible in the following sense:
\begin{eqnarray*}
_M\Delta(hm)=h^{(1)}m^{(-1)}S^{-1}(h^{(3)})\ot h^{(2)}m^{(0)},\\
\mbox{\em if $M$ is a left module and a left comodule ;}\\
 \Delta_M(hm)=h^{(2)}m^{(0)}\ot h^{(3)}m^{(1)}S(h^{(1)}),\\
\mbox{\em if $M$ is a left module and a right comodule ;}\\
_M\Delta(mh)=S(h^{(3)})m^{(-1)}h^{(1)}\ot m^{(0)}h^{(2)},\\
\mbox{\em if $M$ is a right module and a left comodule ;}\\
\Delta_M(mh)=m^{(0)}h^{(2)}\ot S^{-1}(h^{(1)})m^{(1)}h^{(3)},\\
\mbox{\em if $M$ is a right module and a right comodule.}
\end{eqnarray*}
In the first case we say $M$ is stable if $m^{(-1)}m^{(0)}=m$ for all $m
\in M$ (similar definitions apply in other cases).
\end{definition}\noindent

 Let $M$ be an SAYD H-module of the second type in the above definition (left module and right comodule).
 Let $C_n^{alg}(H, M):= M\otimes H^{\otimes (n+1)}$. It is shown in \cite{hkrs} that the
 following operators define a paracyclic
 module structure on $\{ C_n^{alg}(H, M)\}_n$:
 \begin{eqnarray*}
 \delta_i (m\otimes h_0 \otimes \dots\otimes h_n)&=& m\otimes h_0 \otimes\dots
 \otimes h_i h_{i+1} \otimes \dots\otimes h_n,\\
 \delta_n (m \otimes h_0 \otimes\dots\otimes h_n)&=& h_n^{(1)}m \otimes h_n^{(2)}h_0 \otimes \dots\otimes h_{n-1},\\
\sigma_i (m \otimes h_0 \otimes\dots\otimes h_n)&=& m\otimes h_0
\otimes\dots\otimes h_i \otimes 1\otimes h_{i+1}
\otimes \dots h_n,\\
\tau (m \otimes h_0 \otimes\dots\otimes h_n)&=&h_n^{(1)}m \otimes
h_n^{(2)} \otimes h_0 \otimes\dots\otimes h_{n-1}.
 \end{eqnarray*}

 It is shown in \cite{hkrs} that the above operators restrict to the
 subcomplex $C_n^H (H, M):= M\Box_H H^{\otimes (n+1)}$ of  invariant chains on $H$
 with coefficients in $M$ and define a cyclic module that we denote it by $C^H_*(H, M)$.
  Here $\Box $ denotes the cotensor product. We recall that, in general, the cotensor
  product $M\Box_H N$ of a right $H$-comodule $M$ and a left $H$-comodule $N$ is defined
  as the kernel of the map
  $$\Delta_M \otimes 1-1\otimes  _N\Delta :M\otimes N \longrightarrow M\otimes H\otimes N.$$
  The cyclic homology of $H$
  with coefficients in $M$ is by definition the cyclic homology of this
  module. For $M=k$, we are reduced to  the cyclic module
  $\widetilde{H}_{\natural}^{(\delta , \sigma)}$ defined in \cite{ta} and, independently, \cite{kr}.

 With $H$ and $M$ as above, let $C^n_{coalg} (H, M): = H^{\otimes (n+1)} \otimes M$. Endowed with the
 following operators, $\{ C^n_{coalg} (H, M) \}^n$ is a paracocyclic
 module \cite{hkrs}:
 \begin{eqnarray*}
 d_i ( h_0 \otimes \dots\otimes h_n \otimes m)&=&h_0\otimes\dots\otimes
 h_i^{(1)}\otimes h_i^{(2)} \otimes h_{i+1} \otimes\dots\otimes h_n
 \otimes m,\\
 d_{n+1}( h_0 \otimes \dots\otimes h_n \otimes m)&=&h_0^{(2)}\otimes h_1
 \otimes \dots\otimes h_n \otimes h_0^{(1)}S^{-1} (m^{(1)})\otimes
 m^{(2)},\\
 s_i ( h_0 \otimes \dots\otimes h_n \otimes m)&=&h_0 \otimes \dots\otimes
 \epsilon (h_i) \otimes \dots\otimes h_n \otimes m,\\
 t_n (h_0 \otimes \dots\otimes h_n \otimes m)&=&h_1 \otimes h_2\dots\otimes h_n \otimes h_0S^{-1} (m^{(1)})\otimes
 m^{(0)}.
\end{eqnarray*}
To define the cyclic cohomology of $H$ with coefficients in $M$ we
consider the quotient complex $ C^n_H(H, M):= H^{\otimes (n+1)} \otimes_H M$ of invariant cochains
on $H$ with coefficients in $M$. It
is shown in \cite{hkrs} that $ C^*_H(H, M) $ is in fact a cocyclic module. For $M=k$,
we obtain the Connes-Moscovici
cocyclic module $H_{(\delta,\sigma)}^\natural$.

A quick look at the above two modules show that there is no apparent
cyclic duality relationship between $ C^*_H(H, M) $ and $ C^H_*(H, M) $. In the following,
however, we will show,
via a very non-trivial map,   that they are indeed cyclic dual of each other. Consider the dual
paracyclic module $K_* (H, M):= \widehat{C}_{coalg}^* (H, M)$. Using our formulas in
Section 2, its simplicial and cyclic operators are given by:
\begin{eqnarray*}
\delta_i (h_0 \otimes h_1 \otimes \cdots \ot h_n \otimes m)&=& h_0 \otimes
\cdots
\otimes \epsilon (h_i) \otimes  \cdots \ot h_n \otimes m,\\
\sigma_i (h_0 \otimes h_1 \otimes \cdots \ot h_n \otimes m)&=&h_0 \otimes
h_i^{(1)}\otimes h_i^{(2)} \otimes \cdots \otimes h_n \otimes m,\\
\tau (h_0 \otimes h_1 \otimes \cdots \ot h_n \otimes m)&=&h_nm^{(1)}\otimes h_o
\otimes \cdots \otimes h_{n-1} \otimes m^{(0)}.
\end{eqnarray*}

\begin{prop}
The following map defines a morphism of paracyclic modules $\theta : K_* (H,
M) \longrightarrow C_*^{alg} (H, M)$,
\begin{eqnarray*}
 \theta (h_0 \otimes h_1 \otimes \cdots \otimes h_n \otimes m)&=& h_n^{(2)} m^{(0)}
\otimes h_n^{(3)} m^{(1)} S(h_0^{(1)}) \otimes h_0^{(2)} S(h_1
^{(1)})\otimes h_1^{(2)} S (h_2^{(1)}) \otimes  \\
& & \cdots \otimes h_{n-1}^{(2)}
S(h_n^{(1)}).
\end{eqnarray*}
\end{prop}
\begin{proof}
To prove that $\theta$ is a cyclic map one needs to show  that $\theta\delta_i=\delta_i\theta$,
$\theta\sigma_i=\sigma_i\theta$ for $1\le i\le n$, and
$\theta\tau=\tau\theta$. Here  we just check the latter identity  and leave the rest to the reader.
\begin{align*}
&\theta\tau(h_0\ot h_1\ot\dots \ot h_n\ot m)=\theta(h_nm^{(1)}\otimes h_o
\otimes \cdots \otimes h_{n-1} \otimes m^{(0)})\\
&=h_{n-2}\sw{2}m\sw{0}\ot h_{n-1}\sw{3}m\sw{1}S(h_n\sw{1}m\sw{2})\ot h\sw{2}m
\sw{3}S(h_0\sw{1})\ot\dots \ot h_{n-2}\sw{2}S(h_{n-1}\sw{1})\\
&=h_{n-2}^{(2)}m^{(0)}\ot h_{n-1}^{(3)}S(h_n\sw{1})\ot h_n\sw{2}m\sw{1}
S(h_0\sw{1})\ot \dots \ot h_{n-2}\sw{2}S(h_{n-1}\sw{1})\\
&= h_{n-2}^{(2)}S(h_n\sw{2})h_n\sw{3}m^{(0)}\ot h_{n-1}^{(3)}
S(h_n\sw{1})\ot h_n\sw{4}m\sw{1}S(h_0\sw{1})\ot \dots \ot h_{n-2}
\sw{2}S(h_{n-1}\sw{1})\\
&=(h_{n-2}^{(2)}S(h_n\sw{1})\sw{1}h_n\sw{3}m^{(0)}\ot
(h_{n-1}^{(2)}S(h_n\sw{1}))\sw{2}\ot h_n\sw{3}m\sw{1}S(h_0\sw{1})\ot
\dots \ot h_{n-2}\sw{2}S(h_{n-1}\sw{1})\\
&=\tau\theta(h_0\ot h_1\ot\dots \ot h_n\ot m).
\end{align*}

\end{proof}

The map $\theta $ is neither injective nor surjective. We can however
show that it descends to the relevant spaces of invariants on both sides and gives an isomorphism of cyclic modules:

\begin{lem} The map $\theta$ lands in $M\Box_H H^{\ot(n+1)}$ and  descends to the
quotients $H^{\otimes (n+1)} \otimes_H M$ for each $n $.
\end{lem}
\begin{proof}
To prove  the second part  we show that for all  $g, h_0,h_1,\dots h_n\in H$ and $m\in M$ we have
$\theta((h_0\ot h_1\ot\dots\ot h_n)\cdot g\ot m)=\theta(h_0\ot h_1\ot\dots\ot h_n\ot gm)$.
 Indeed:
\begin{align*}
  &\theta((h_0\ot h_1\ot\dots\ot h_n)\cdot g\ot m)=\theta(h_0g\sw{1}\ot h_1g\sw{2}\ot\dots\ot h_ng\sw{n+1}\ot m)\\
  &=(h_ng\sw{n+1})^{(2)} m^{(0)}
\otimes (h_ng\sw{n+1})^{(3)} m^{(1)} S((h_0g\sw{1})\sw{1}) \otimes (h_0g\sw{1})^{(2)} S((h_1g\sw{2})
^{(1)})\otimes\\
& (h_1g\sw{2})^{(2)} S ((h_2g\sw{3})^{(1)}) \otimes  \cdots \otimes (h_{n-1}g\sw{n})^{(2)}S((h_ng\sw{n+1})^{(1)})\\
&=h_n\sw{2}g\sw{2n+2} m^{(0)}
\otimes h_n\sw{3}g\sw{2n+3} m^{(1)} S(g\sw{1})S(h_0\sw{1}) \otimes h_0\sw{2}g\sw{2} S(g\sw{3})S(h_1\sw{1})
\otimes \\
&h_1\sw{2}g\sw{4} S (g\sw{5})S(h_2\sw{1}) \otimes   \dots \otimes (h_{n-1}\sw{2}g\sw{2n}S(g\sw{2n+1})S(h_n\sw{1})\\
&=h_n\sw{2}g\sw{2}m\sw{0}\ot h_n\sw{3}g\sw{3}m\sw{1}S(g\sw{1})\ot h_0\sw{2}S(h_1\sw{1})\ot
\dots\ot h_{n-1}\sw{2}S(h_n\sw{1}).
\end{align*}
On the other hand, by using the stable anti-Yetter-Drinfeld module property of $M$, one has:
\begin{align*}
&\theta(h_0\ot h_1\ot\dots\ot h_n\ot gm)=\\
&h_n\sw{2}(gm)\sw{0}\ot h_n\sw{3}(gm)\sw{1}S(h_0\sw{1})\ot h_0\sw{2}S(h_1\sw{1})\ot
\dots\ot h_{n-1}\sw{2}S(h_n\sw{1})\\
&h_n\sw{2}g\sw{2}m\sw{0}\ot h_n\sw{3}g\sw{3}m\sw{1}S(g\sw{1})\ot h_0\sw{2}S(h_1\sw{1})\ot \dots\ot
h_{n-1}\sw{2}S(h_n\sw{1}).
\end{align*}

Now we prove that the image of $\theta$ is in  the cotensor product space  $M\Box_HH^{\ot(n+1)}$.
That is,
$(\Delta_M\ot id_{H^{\ot(n+1)}})\circ\theta=(id_M\ot _{H^{\ot(n+1)}}\hspace{-2pt}\Delta)\circ\theta$.
\begin{align*}
&\Delta_M\ot id_{H^{\ot(n+1)}}\circ\theta(h_0\ot h_1\ot\dots\ot h_n\ot m)\\
&=\Delta_M\ot id_{H^{\ot(n+1)}}( h_n^{(2)} m^{(0)}
\otimes h_n^{(3)} m^{(1)} S(h_0^{(1)}) \otimes h_0^{(2)} S(h_1
^{(1)})\otimes h_1^{(2)} S (h_2^{(1)}) \otimes
 \cdots \\
 &\otimes h_{n-1}^{(2)}S(h_n^{(1)}))\\
 &=(h_n^{(2)} m^{(0)})\sw{0}\ot (h_n^{(2)} m^{(0)})\sw{1}\otimes h_n^{(3)} m^{(1)} S(h_0^{(1)})
 \otimes h_0^{(2)} S(h_1
^{(1)})\otimes h_1^{(2)} S (h_2^{(1)}) \otimes
 \cdots\\
 &\otimes h_{n-1}^{(2)}S(h_n^{(1)}))\\
 &=h_n^{(3)} m^{(0)}
\otimes h_n^{(4)} m^{(1)} S(h_n^{(2)})\ot h_n\sw{5}m\sw{2}S(h_0\sw{1}) \otimes h_0^{(2)} S(h_1
^{(1)})\otimes  \cdots \otimes h_{n-1}^{(2)}S(h_n^{(1)}))
\end{align*}
On the other hand, we have:
\begin{align*}
&id_M\ot _{H^{\ot(n+1)}}\hspace{-2pt}\Delta\circ\theta(h_0\ot h_1\ot\dots\ot h_n\ot m)=\\
&=h_n\sw{2}m\sw{0}\ot(h_n\sw{3}m\sw{1}S(h_0\sw{1}))\sw{1}(h_0\sw{2}S(h_1\sw{1}))\sw{1}
\dots(h_{n-1}\sw{2}S(h_n\sw{1}))\sw{1}\ot\\
&(h_n\sw{3}m\sw{1}S(h_0\sw{1}))\sw{2}\ot(h_0\sw{2}S(h_1\sw{1}))\sw{2}\ot\dots\ot
(h_{n-1}\sw{2}S(h_n\sw{1}))\sw{2}\\
&=h_n^{(3)} m^{(0)}
\otimes h_n^{(4)} m^{(1)} S(h_n^{(2)})\ot h_n\sw{5}m\sw{2}S(h_0\sw{1}) \otimes h_0^{(2)} S(h_1
^{(1)})\otimes  \cdots \otimes h_{n-1}^{(2)}S(h_n^{(1)})).
\end{align*}
\end{proof}

\begin{theorem} The induced map $\theta$ is an isomorphism of cyclic
modules
$$ \theta : K_*^H (H, M) \longrightarrow C_*^H(H, M).$$
\end{theorem}
\begin{proof}
At  first we identify  $H^{\ot(n+1)}\ot_HM$ with $H^{\ot n}\ot M$.
Consider the map $\phi: H^{\ot(n+1)}\otimes M\rightarrow H^{\ot n}\ot
M$,
$$\phi(h_0\ot h_1\ot\dots\ot h_n\ot m)=h_0S(h_n\sw{n})\ot h_1S(h_n\sw{2})
\ot\dots\ot h_{n-1}S(h_n\sw{1})\ot h_n\sw{n+1} m.$$
One can see that this map is $H$-bilinear. So it induces a  map
$\bar\phi :H^{\ot(n+1)}\ot_HM \rightarrow H^{\ot n}\ot M$.
 It can easily be checked that the following map is the two-sided
 inverse of $\bar\phi$:
\begin{align*}
&\psi: H^{\ot n}\ot M\rightarrow H^{\ot (n+1)}\ot_H M,\\
&\psi(h_1\ot h_2\ot\dots\ot h_n\ot m)= h_1\ot h_2\ot\dots\ot h_n\ot 1\ot_H m.
 \end{align*}

  Next we identify $M\Box_HH^{\ot(n+1)}$ with $M\ot H^{\ot n}$. Consider the  maps
  $\phi':M\Box_HH^{\ot(n+1)}\longrightarrow M\ot H^{\ot n}$,
 $$\phi'(m\ot h_0\ot \dots\ot h_n)=m\ot\epsilon(h_0)h_1\ot h_2\ot
 \dots\ot h_n,$$
  and $\psi':M\ot H^{\ot n}\longrightarrow M\Box_HH^{\ot(n+1)}$,
  $$\psi'(m\ot h_1\ot h_2\ot\dots \ot h_n)=
  m\sw{0}\ot m\sw{1}S(h_1\sw{1}h_2\sw{1}\dots h_n\sw{1})\ot h_1\sw{1}\ot h_2\sw{1}\ot
  \dots\ot  h_n\sw{1}.$$
   One can check that $\phi'$ and $\psi'$ are inverse to one another.

 After  the above identifications,  the map induced  by $\theta$,
 denoted  $\bar\theta$, on $H^{\ot(n+1)}\ot_HM$, has the following
  formula:
  \begin{align*}
  &\bar\theta:H^{\ot n}\ot M\longrightarrow M\ot H^{\ot n},\\
  &\bar\theta(h_1\ot h_2\ot\dots \ot h_n\ot m )=m\sw{0}\ot m\sw{1}S(h_0\sw{1})\ot
  h_0\sw{2}S(h_1\sw{1})\ot\dots \\
  &\ot h_{n-2}\sw{2}S(h_{n-1}\sw{1})\ot h_{n-1}\sw{2}.
  \end{align*}
  One can directly check  that $\bar\theta $ is an isomorphism and its inverse is
  given by
  \begin{align*}
  &\gamma:  M\ot H^{\ot n}\longrightarrow H^{\ot n}\ot M,\\
  &\gamma(m\ot h_1\ot h_2\ot\dots \ot h_n)= S(m\sw{1})h_0h_1\sw{1}\dots h_{n-1}\sw{1}\ot h_1\sw{2}\dots
  h_{n-1}\sw{2}\ot
  \dots\\
  &\ot h_{n-2}\sw{n-1}h_{n_2}\sw{n-1}\ot h_{n-1}\sw{n}\ot m^{(0)}.
  \end{align*}

\end{proof}
The cyclic module $K_*^H(H, M)$ is used by Jara and Stefan  in their
study of relative cyclic homology of Hopf-Galois extensions \cite{js}. Note that what we call an stable
anti-Yetter-Drinfeld module in the present paper (and in \cite{hkrs}), is called a modular crossed module in
\cite{js}. It follows
from Theorem 3.1 above that Theorem 4.13 in \cite{js} is a consequence
of Theorem 3.1 in \cite{hkrs} (by choosing $A=H$). For the same reason, Theorem 5.2 in \cite{js}
follows from Theorem 3.22
in our paper \cite{kr5} (for special classes of SAYD modules called matched pairs, but the same
 proof works in general).

In the remainder of this paper we briefly look into a different type of duality, i.e. the
$Hom$-functor duality,
 between the  cyclic module $C^H_*(H, M)$ and the cocyclic module
 $C^*_H(G, N)$. Recall that a
 Hopf pairing between Hopf algebras  $G$ and $H$ is a
 bilinear map
$<, >: H\otimes G \longrightarrow k$ that satisfies the following
relations for all $h, h_1, h_2$ in $H$ and $g, g_1, g_2$ in $G$:
$$ <h_1 h_2, g>= <h_1, g^{(1)}> <h_2, g^{(2)}>, \quad <h, g_1 g_2>=<h^{(1)}, g_1><h^{(2)}, g_2>,$$
$$<h,1>= \varepsilon (h), \quad \quad <1, g>= \varepsilon (g).$$
 In addition let $M\in\ _H\hspace{-2pt}\mathcal{M}^H$ be a left $H$-moodule and a right $H$-comodule, and $
 N\in\ ^G\hspace{-2pt}\mathcal{M}_G$ be a left $G$-comodule and a right $G$-module, both
 satisfying the SAYD-module conditions.
 To complete the set up for having
 a pairing between  Hopf cyclic cohomologies we need a
 pairing between $M$ and $N$, denoted here by $ <,>: M\otimes N \rightarrow k$, such that for all
  $m\in M,\; n\in N,\;h\in H,\;g\in  G$, the following
 conditions are satisfied:
 $$ <hm,n>=<h,n\sw{-1}><m,n\sw{0}>,\quad\quad <m,ng>=<m\sw{0},n><m\sw{1},g>.$$
Consider the map $ M\otimes H^{\otimes n} \longrightarrow Hom (N\otimes G^{\otimes
n}, k)$ defined by
$$ (m\ot h_1 \otimes \dots \otimes h_n) (n\ot g_1\otimes  \dots \otimes g_n)=
<m,n>\prod_{i=1}^{i=n}<h_i , g_i>.$$

\begin{prop} The above map defines  a morphism of cocyclic modules
$$C^\ast_{coalg}(H,M) \longrightarrow Hom_k(C_\ast^{alg}(G,N),k).$$
\end{prop}

 If $H$ is finite dimensional over a field $k$ and $G=H^*$ is the dual Hopf algebra, then it is
 easy to see that the above map, for $M=N=k$ and the natural pairing
 between $H$ and $H^*$,
 is an isomorphism of cocyclic modules. Thus the Connes-Moscovici
 cyclic cohomology of a finite dimensional Hopf algebra
 is isomorphic
 to the cyclic cohomology in the sense
 of \cite{kr, ta} of  $G=H^*$. In the infinite dimensional case, however,  the induced
  map on cohomology can be trivial, even
 for a non-degenerate pairing. For example, let $H=Rep( {\bf G})$ be the
 Hopf algebra of representable functions on a compact Lie group ${\bf G}$ and $G=
 U(\mathfrak{g})$ the enveloping algebra of the Lie algebra $\mathfrak{g}$ of
 ${\bf G}$. The canonical pairing $H \otimes G \rightarrow \mathbb{C}$
  defined by
 $$<f, X_1\otimes \dots \otimes X_n >= \frac{d}{dt} f(e^{tX_1}\dots e^{tX_n})\vert_{t=0}$$
 is non-degenerate. Then while the Connes-Moscovici cyclic cohomology of $H$ is
 trivial, the cyclic cohomology of $G$ in the sense of \cite{kr, ta} is non-trivial and in fact is
 isomorphic to the Lie algebra
 cohomology of $\mathfrak{g}$ (see \cite{kr} for a proof of both statements).

  In Lemma 2.3 we saw that cyclic duals of cyclic modules of algebras
  and coalgebras are both homologically trivial. It follows from Theorem 3.1 that this need not
  be true for cyclic modules of
  Hopf algebras.  In fact for $H=U(\mathfrak{g})$ and the modular pair $(\delta , \sigma )= (\varepsilon , 1)$,
   the periodic cyclic cohomology and homology are computed in \cite{achm98} and \cite{kr}
   respectively. They are both
    isomorphic to the Lie algebra homology of $\mathfrak{g}$.

\end{document}